\documentclass[draft|showlink]{siamart1116}

\usepackage{amsfonts}
\usepackage{graphicx}
\usepackage{epstopdf}
\usepackage{algorithmic}
\ifpdf
 \DeclareGraphicsExtensions{.eps,.pdf,.png,.jpg}
\else
 \DeclareGraphicsExtensions{.eps}
\fi

\newtheorem{scheme}{Scheme}
\newtheorem{remark}{Remark}
\newtheorem{example}{Example}

\newcommand{\TheTitleABB}{An efficient third-order scheme for BSDEs} 
\newcommand{\TheTitle}{An efficient third-order scheme for BSDEs based on nonequidistant difference scheme} 
\newcommand{\TheAuthors}{Chol-Kyu Pak, Mun-Chol Kim and Chang-Ho Rim }

\def\R{\mathbb R} \def\P{\mathbb P} \def\F{\mathcal F} \def\d{\partial} \def\to{\rightarrow} \def\x{\mathbf x}
\def\beq{\begin{equation}} \def\enq{\end{equation}}
\def\beseq{\begin{subequations}} \def\enseq{\end{subequations}}
\def\beqa{\begin{eqnarray}} \def\enqa{\end{eqnarray}}

\def\BeDef{\begin{definition}} \def\EnDef{\end{definition}}
\def\BeThe{\begin{theorem}} \def\EnThe{\end{theorem}}
\def\BeLem{\begin{lem}} \def\EnLem{\end{lem}}

\numberwithin{theorem}{section}
\numberwithin{equation}{section}
\numberwithin{table}{section}
\numberwithin{figure}{section}
\numberwithin{remark}{section}
\numberwithin{example}{section}

\headers{\TheTitleABB}{\TheAuthors}

\title{{\TheTitle}}

\author{
 Chol-Kyu Pak\thanks{Faculty of Mathematics, Kim Il Sung University, Pyongyang, Democratic People's Republic of Korea (\email{pck2016217@gmail.com}).}
\and
 Mun-Chol Kim
\and
Chang-Ho Rim
}

\usepackage{amsopn}

\begin{document}

\maketitle

% REQUIRED
\begin{abstract}
In this paper we propose an efficient third-order numerical scheme for backward stochastic differential equations(BSDEs). We use 3-point Gauss-Hermite quadrature rule for approximation of the conditional expectation and avoid spatial interpolation by setting up a fully nested spatial grid and using the approximation of derivatives based on non-equidistant sample points. As a result, the overall computational complexity is reduced significantly. Several examples show that the proposed scheme is of third-order and very efficient.
\end{abstract}

% REQUIRED

% REQUIRED

\section{Introduction}\label{sec:sec1}
Let $(\Omega,\F, \P)$ be a probability space, $T>0$ a finite time and $\{\F_t\}_{0\le t \le T}$ a filtration satisfying the usual conditions. Let $(\Omega,\F, \P,\{\F_t\}_{0\le t\le T})$ be a complete filtered probability space on which a standard $d$-dimensional Brownian motion $W_t=(W_t^1,W_t^2, \cdots ,W_t^d )^T$ is defined and $\F_0$ contains all the $\P$-null sets of $\F$. 
\par The general form of backward stochastic differential equation (BSDE) is
\beq\label{eq:bsde}
y_t=\xi+\int_t^T{f(s,y_s,z_s )ds}-\int_t^T{z_sdWs}, \quad t\in[0,T] 
\enq
where the generator $f:[0,T]\times\R^m\times\R^{m\times d}\to \R^m$ is $\{\F_t\}$-adapted for each $(y,z)$ and the terminal variable $\xi$ is a $\F_T$-measurable and square integrable random variable.
A process $(y_t,z_t ):[0,T]\times\Omega \rightarrow \R^m\times \R^{m\times d}$ is called an $L^2$-solution of the BSDE {\cref{eq:bsde} if it satisfies the equation \cref{eq:bsde} while it is ${\{\F_t\}}$-adapted and square integrable .
\par In 1990, Pardoux and Peng first proved in \cite{Peng90} the existence and uniqueness of the solution of general nonlinear BSDEs and afterwards there has been very active research in this field with many applications.(\cite{Karoui97})
\par In this paper we assume that the terminal condition is a function of $W_T$, i.e., $\xi=\varphi(W_T)$ and the BSDE \cref{eq:bsde} has a unique solution $(y_t,z_t )$. 
It was shown in \cite{Peng91} that the solution $(y_t,z_t)$ of \cref{eq:bsde} can be represented as
\beq\label{eq:eq12}
y_t=u(t,W_t ), z_t=\nabla_x u(t,W_t ), \; \forall t\in[0,T)
\enq
where $u(t,x)$ is the solution of the parabolic partial differential equation
\beq\label{eq:pde}
\frac{\d u}{\d t}+\frac{1}{2} \sum_{i=1}^d \frac{\d^2 u}{\d x_i^2}+f(t,u,\nabla_x u)=0
\enq
with the terminal condition $u(T,x)=\varphi (x)$, and $\nabla_x u$ is the gradient of $u$ with respect to the spatial variable $x$. The smoothness of $u$ depends on $f$ and $\varphi$.
\par Although BSDEs and their extensions such as FBSDEs have very important applications in many fields such as mathematical finance and stochastic control, it is well known that it is difficult to obtain the analytical solutions except some special cases and there have been many works on numerical methods.(\cite{Bou, Gobet, Prot02, Prot94, Tre06, Zhao06, Ruijter, Zhao14, Zhao17}) 
\par Among all the previous works, we are concerned with \cite{Zhao14} where a new kind of multi-step scheme for FBSDEs was proposed. They demonstrated experimentally that the scheme is of high order (up to 6th order) and that scheme was simpler than the previous ones such as \cite{Zhao06, Zhao10}. 
\par Our scheme can be seen as a continuation of \cite{Zhao14} for standard BSDEs driven by a Wiener process, but our scheme includes some important ideas. 
\par The approximation of conditional expectation is a common problem in solving BSDEs numerically and it influences both the convergence order and the computational complexity. Gauss-Hermite quadrature is widely used for this purpose because it can achieve high accuracy using only a few sample points. At the same time one needs the approximation values at non-grid points and some kinds of interpolation is required. This contributes much to the computational complexity. So some efficient methods such as sparse-grid method has been developed.(\cite{Sparse})
\par In this paper we propose a new efficient third-order scheme for BSDEs which does not imply interpolations. The main idea of new scheme is to make the time-space grid include all quadrature points by using 3-point Gauss-Hermite quadrature rule and the derivative approximation based on non-equidistant sample points. 
As a result, the computational cost can be reduced significantly because the spatial interpolation is not needed at all while third-order convergence is achieved, although it still suffers from curse of dimensionality for multidimensional problems.
\par We note that the previous works used a big number(e.g. 8 or 10) for the number of Gauss-Hermite quadrature points for each dimension and a high-order polynomial interpolation(e.g. 12 or 15th order).(\cite{Zhao10, Zhao14, Zhao17}) It contributed much to computational complexity and our estimates show how to set the number of Gauss-Hermite quadrature points  properly.
\par The nonquidistant difference scheme for derivative approximation is very stable for high-order and this kind of difference scheme could have applications in other fields.
\par The rest of the paper is organized as follows. In \Cref{sec:sec2}, we derive a discrete scheme for BSDEs using the idea in \cite{Zhao14}. In \Cref{sec:sec3}, we propose a new efficient third-order scheme based on nonequidistant difference scheme through some analysis on the approximation of conditional expectation by Gauss-Hermite quadrature. In \Cref{sec:sec4}, we demonstrate the convergence rate and efficiency of our new scheme through several examples and finally some conclusions are given in \Cref{sec:sec5}.
\section{Discrete scheme based on derivative approximation}\label{sec:sec2}
Let us consider the following backward stochastic differential equation (BSDE):
\beq\label{eq:eq21} y_t=\varphi(W_T)+\int_t^T{f(s,y_s,z_s)ds}-\int_t^T{z_sdW_s}, \quad t\in [0,T] \enq
where the generator $f:[0,T]\times\Omega\times\R^m\times\R^{m\times d}\to\R^m$ is a stochastic process that is $\{\F_t\}$-adapted for all $(y,z)$ and $\varphi:\R^d\to\R^m$ is a measurable function. As in \cite{Zhao10} we assume that $f$ and $\varphi$ are smooth enough and their derivatives as well as themselves are all bounded.
\par Let $0=t_0<\cdots<t_N=T$ be an equidistant partition of $[0,T]$ and $t_{n+1}-t_n=h=T/N$.
As in \cite{Zhao06,Zhao10}, we define $\F_s^{t,\mathbf{x}}(t\leq s\leq T)$ to be be a $\sigma$-field generated by the Brownian motion $\{\mathbf{x}+W_r-W_t,t\leq r\leq s\}$ starting from the time-space point $(t,\mathbf{x})$ and $E_t^\mathbf{x}[\cdot]:=E[\cdot|\F_t^{t,\mathbf{x}}]$.
Taking the conditional expectation $E_{t_n}^\mathbf{x}[\cdot]$ on the both sides of \cref{eq:eq21}, we get
\beq\label{eq:eq22} E_{t_n}^\mathbf{x}[y_t]=E_{t_n}^\mathbf{x}[\varphi(W_T)]+\int_t^T{E_{t_n}^\mathbf{x}[f(s,y_s,z_s)]ds}\enq
By differentiating the both sides of \cref{eq:eq22} with respect to $t$, we obtain the following ODE.(see \cite{Zhao14} for differentiability)
\beq\label{eq:eq23} \frac{dE_{t_n}^\mathbf{x}[y_t]}{dt}=-E_{t_n}^\mathbf{x}[f(s,y_s,z_s)]\enq
Next, for $t\in[t_n,T]$ we have
\beq\label{eq:eq24} y_{t_n}=y_t+\int_t^T{f(s,y_s,z_s)ds}-\int_{t_n}^t{z_sdW_s}.\enq
Multiplying $\Delta W_{t_n,t}^\mathbf{T}:=W_t^\mathbf{T}-W_{t_n}^\mathbf{T}$ to the both sides of \cref{eq:eq24} and taking $E_{t_n}^\mathbf{x}[\cdot]$, we obtain
\beq\label{eq:eq25} 0=E_{t_n}^\mathbf{x}[y_t\Delta W_{t_n,t}^\mathbf{T}]+\int_{t_n}^t{E_{t_n}^\mathbf{x}[f(s,y_s,z_s)\Delta W_{t_n,t}^\mathbf{T}]ds}-\int_{t_n}^t{E_{t_n}^\mathbf{x}[z_s]ds}\enq
where $(\cdot)^\mathbf{T}$ means transpose of $(\cdot)$.
Taking derivatives of the both sides of \cref{eq:eq25} with respect to $t$, we obtain the following ODE.
\beq\label{eq:eq26} E_{t_n}^\mathbf{x}[z_t]=\frac{dE_{t_n}^\mathbf{x}[y_t\Delta W_{t_n,t}^\mathbf{T}]}{dt}+E_{t_n}^\mathbf{x}[f(t,y_t,z_t)\Delta W_{t_n,t}^\mathbf{T}]\enq
The two ODEs \cref{eq:eq23} and \cref{eq:eq26} are called the reference equations for the BSDE \cref{eq:eq21}.
From these reference equations we have
\beq\label{eq:eq27} \frac{dE_{t_n}^\mathbf{x}[y_t]}{dt}|_{t=t_n}=-f(t_n,y_{t_n},z_{t_n})\enq
\beq\label{eq:eq28} z_{t_n}=\frac{dE_{t_n}^\mathbf{x}[y_t\Delta W_{t_n,t}^\mathbf{T}]}{dt}|_{t=t_n}.\enq
Now we approximate the derivatives in \cref{eq:eq27} and \cref{eq:eq28}. In \cite{Zhao14} the approximation formula was derived by solving linear system obtained by Taylor's expansion. In this paper we approximate the derivatives by the ones of Lagrange interpolaion polynomial based on equidistant sample points and still get the same result.
\par Let $u(t):\R\to\R$ be $k+1$ times differentiable, then the Lagrange interpolation polynomial based on values $\{u_0,u_1,\cdots,u_k\}$ on equidistant sample points $\{t_0,t_0+h,\cdots,t_0+kh\}$ can be written as
\beq\label{eq:eq29} L(t)=\sum_{i=0}^k{\frac{\prod_{i\neq j}{(t-t_0-jh)}}{\prod_{i\neq j}{(i-j)}}h^{-k}u_i}\enq
and the deviation is given by
\beq\label{eq:eq210} L(t)-u(t)=\frac{f^{(k+1)}(\xi)\prod_{i=0}^n{(t-t_0-ih)}}{(k+1)!}, \quad t_0\leq \xi \leq t_n.\enq
By differentiating \cref{eq:eq29}, we get
\beq\label{eq:eq211} \frac{dL}{dt}(t)=\sum_{i=0}^k{\frac{h^{-k}u_i}{\prod_{j\neq i}{(i-j)}}\sum_{j\neq i}{\prod_{l\neq i, l\neq j}{(t-t_0-lh)}}}\enq
and furthermore,
\begin{multline}\label{eq:eq212}
\begin{aligned}
\frac{dL}{dt}(t_0)&=\sum_{i=0}^k{\frac{(-1)^{k-1}h^{-1}u_i}{\prod_{j\neq i}{(i-j)}}\sum_{j\neq i}{\prod_{l\neq i, l\neq j}{l}}}\\&=h^{-1}\left(-\sum_{j\neq0}{\frac{1}{j}}u_0+\sum_{i=1}^k{\frac{(-1)^{k-1}k!}{i\prod_{i\neq j}{(i-j)}}}u_i\right)\\&=h^{-1}\left(-\sum_{j\neq0}{\frac{1}{j}}u_0+\sum_{i=1}^k{\frac{(-1)^{i-1}\begin{pmatrix}k\\i\end{pmatrix}}{i}}u_i\right).
\end{aligned}
\end{multline} 
Let $\alpha_i^k$ be the coefficients of $u_i$ in \cref{eq:eq212}, then we have
\beq\label{eq:eq213} \frac{dL}{dt}(t_0)=\frac{\sum_{i=0}^k{\alpha_i^k u_i}}{h}.\enq
In \Cref{tbl21} we list $\{\alpha_i^k\}$ for $k=1\cdots6$.
\begin{table}[tbhp]\label{tbl21}
\caption{$\{\alpha_i^k\}$ in the derivative approximation based on equidistant Lagrange interpolation}
\centering
\begin{tabular}{|c|c|c|c|c|c|c|c|} \hline
$\alpha_i^k$& $i=0$ & $i=1$ & $i=2$ & $i=3$ & $i=4$ & $i=5$ & $i=6$ \\ \hline
$k=1$ &-1&1&&&&&\\ \hline
$k=2$ & $-\frac{3}{2}$ &2& $\frac{1}{2}$ &&&&\\ \hline
$k=3$ & $-\frac{11}{6}$ &3& $-\frac{3}{2}$ & $\frac{1}{3}$ &&&\\ \hline
$k=4$ & $-\frac{25}{12}$ &4&-3& $\frac{4}{3}$ & $-\frac{1}{4}$ &&\\ \hline
$k=5$ & $-\frac{137}{60}$ &5&-5& $\frac{10}{3}$ & $-\frac{5}{4}$ & $\frac{1}{5}$ &\\ \hline
$k=6$ & $-\frac{49}{20}$ &6& $-\frac{15}{2}$ & $\frac{20}{3}$ & $-\frac{15}{4}$ & $\frac{6}{5}$ & $-\frac{1}{6}$ \\ \hline
\end{tabular}
\end{table}
From \cref{eq:eq210}, there exists $\xi\in\R$ such that
\beq\label{eq:eq214} \frac{du}{dt}=\frac{\sum_{i=0}^k{\alpha_i^k u_i}}{h}+\frac{(-1)^ku^{(k+1)}(\xi)}{k}h^k.\enq
Note that the coefficients $\{\alpha_i^k\}$ in \Cref{tbl21} are consistent with \cite{Zhao14}.

Finally by approximating the derivatives in \cref{eq:eq27} and \cref{eq:eq28} using \cref{eq:eq214}, we obtain the semi-discrete scheme for BSDE \cref{eq:eq21} as follows.

\begin{scheme}\label{scheme1}
Assume that $y^{N-i},z^{N-i}(i=0,\cdots,k-1)$ are known. For $n=N-k,\cdots,0$ solve $y^n=y^n(\mathbf{x}),z^n=z^n(\mathbf{x})$ at time-space point $(t_n,\mathbf{x})$ by
\beq\label{eq:eq215} z^n=\frac{\sum_{j=1}^k{\alpha_j^kE_{t_n}^\mathbf{x}[y^{n+j}\Delta W_{t_n,t_{n+j}}^\mathbf{T}]}}{h}\enq
\beq\label{eq:eq216} -\alpha_0^ky^n=\sum_{j=1}^k{\alpha_j^kE_{t_n}^\mathbf{x}[y^{n+j}]}+h\cdot f(t_n,y^n,z^n)\enq
\end{scheme}
The local truncation error of the above scheme is $O(h^k)$ from \cref{eq:eq214}.
The \Cref{scheme1} is called a semi-discrete scheme because it is discretized only in time domain. One needs spatial grids to get the fully discrete scheme. Let $D^n\subset \R^d$ be a spatial grid for $t_n$. (The detailed structure of $D^n$ depends on how to approximate the conditional expectations.) The fully-discrete scheme on time-space grid $\bigcup_{i=0}^N{\{t_i\}\times D^i}$ based on \Cref{scheme1} is as follows.

\begin{scheme}\label{scheme2}
Assume that $y^{N-i},z^{N-i}$ on $D^{N-i}$ are known for $i=0,\cdots,k-1$. For $n=N-k,\cdots,0,\mathbf{x}\in D^n$ solve $y^n=y^n(\mathbf{x}),z^n=z^n(\mathbf{x})$ by
\beq\label{eq:eq217} z^n=\frac{\sum_{j=1}^k{\alpha_j^k\hat{E}_{t_n}^\mathbf{x}[y^{n+j}\Delta W_{t_n,t_{n+j}}^\mathbf{T}]}}{h}\enq
\beq\label{eq:eq218} -\alpha_0^ky^n=\sum_{j=1}^k{\alpha_j^k\hat{E}_{t_n}^\mathbf{x}[y^{n+j}]}+h\cdot f(t_n,y^n,z^n)\enq
\end{scheme}

In the above scheme, $\hat{E}_{t_n}^\mathbf{x}[\cdot]$ stands for the approximation of conditional expectation which can be calculated by Monte-Carlo methods or other quadrature rules. We are interested in the case where Gauss-Hermite quadrature rule is used.

\section{Analysis and improvement}\label{sec:sec3}
\par In this section, we discuss how to improve the \Cref{scheme2} in the case where Gauss-Hermite quadrature rule is used for the approximation of conditional expectations. The error of \Cref{scheme2} comes from three sources: the time discretization through derivative approximation, Gauss-Hermite quadrature and spatial interpolation.
\subsection{Error analysis for the approximation of conditional expectation}\label{sec:sec31}
\par We first consider the case of $m=d=1$.
For a function $g:\R\to\R$, the integration $\int_\R{g(x)e^{-x^2}dx}$ can be approximated by Gauss-Hermite quadrature as follows.
\beq\label{eq:eq31} \int_\R{g(x)e^{-x^2}dx}\approx\sum_{i=1}^L{\omega_i g(a_i)}\enq
where $L$ is a parameter of quadrature which stands for the number of sample points used, $\{a_i\}_{i=1}^L$ are the roots of the Hermite polynomial of degree $L$ defined by 
\[ H_L(x)=(-1)^Le^{x^2}\frac{d^L}{dx^L}(e^{-x^2})\]
 , and the weights $\{\omega_i\}_{i=1}^L$ are defined by
\beq\label{eq:eq32} \omega_i=\frac{2^{L+1}L!\sqrt{\pi}}{(H'_L(a_i))^2}.\enq
The truncation error $R(g,L)$ is
\beq\label{eq:eq33} R(g,L)=\int_\R{g(x)e^{-x^2}dx}-\sum_{i=1}^L{\omega_ig(a_i)}=\frac{L!\sqrt\pi}{2^L(2L)!}g^{(2L)}(\bar x),\quad \bar x\in\R.\enq
The previous works such as \cite{Zhao06, Zhao10, Zhao14, Zhao17} paid attention to the coefficient $\frac{L!\sqrt\pi}{2^L(2L)!}$ in the above equality. In fact this coefficient becomes very small for $L\ge8$ (less than about $1.33\times10^{-11}$). So in the above literatures $L=8$ was considered to be satisfactory. (In \cite{Zhao17} they used $L=10$, where the coefficient becomes less than about $2.58\times 10^{-15}$.) \par Here we investigate the error in more detail.
From \cref{eq:eq12}, the solution $(y_t,z_t)$ can be represented as $y_t=u(t,W_t),z_t=\nabla u(t,W_t)$ and we have
\[ E_{t_n}^x[y_{t_{n+j}}]=\frac{1}{\sqrt{2\pi jh}}\int_\R{u(t_{n+j},x+v)e^{-\frac{v^2}{2jh}}dv}=\frac{1}{\sqrt\pi}\int_\R{u(t_{n+j},x+\sqrt{2jh}w)e^{-w^2}dw} \]
where we applied the change of variable $v=\sqrt{2jh}\cdot w$.
Now applying Gauss-Hermite quadrature rule to the kernel $g_1(w)=u(t_{n+j},x+\sqrt{2jh}w)=y_{t_{n+j}}(x+\sqrt{2jh}w_i)$ we have
\beq\label{eq:eq34} E_{t_n}^x[y_{t_{n+j}}]\approx \hat E_{t_n}^x[y_{t_{n+j}}]=\frac{1}{\sqrt{\pi}}\sum_{i=1}^L{\omega_i y_{t_{n+j}}(x+\sqrt{2jh}a_i)}\enq
and the truncation error is 
\begin{multline}\label{eq:eq35}
\begin{aligned} 
\frac{L!\sqrt{\pi}}{2^L(2L)!}g_1^{(2L)}(\bar w)&=\frac{L!\sqrt{\pi}}{2^L(2L)!}\cdot \frac{\d^{2L}u}{\d x^{2L}}(t_{n+j},x+\sqrt{2jh}\bar w)\cdot (2jh)^L\\
&=\frac{j^LL!\sqrt{\pi}}{(2L)!}\cdot \frac{\d^{2L}u}{\d x^{2L}}(t_{n+j},x+\sqrt{2jh}\bar w)\cdot h^L, \quad \bar w\in\R.
\end{aligned}
\end{multline}
Therefore under the assumption that $\frac{\d^{2L}u}{\d x^{2L}}$ is bounded, we obtain
\beq\label{eq:eq36}\left|E_{t_n}^x[y_{t_{n+j}}]-\hat E_{t_n}^x[y_{t_{n+j}}]\right|\leq C(jh)^L.\enq
Similarly we deduce
\begin{multline}\label{eq:eq37}
\begin{aligned} 
E_{t_n}^x[y_{t_{n+j}}\Delta W_{t_n,t_{n+j}}]&=\frac{1}{\sqrt{2\pi jh}}\int_\R{u(t_{n+j},x+v)ve^{-\frac{v^2}{2jh}}dv}\\&=\frac{\sqrt{2jh}}{\sqrt\pi}\int_\R{u(t_{n+j},x+\sqrt{2jh}w)we^{-w^2}dw}
\end{aligned}\end{multline}
and applying Gauss-Hermite quadrature rule we have
\beq\label{eq:eq38} E_{t_n}^x[y_{t_{n+j}}\Delta W_{t_n,t_{n+j}}]\approx \hat E_{t_n}^x[y_{t_{n+j}}\Delta W_{t_n,t_{n+j}}]=\sqrt{\frac{2jh}{\pi}}\sum_{i=1}^L{\omega_i y_{t_{n+j}}(x+\sqrt{2jh}a_i)a_i}\enq
where the truncation error is 
\beq\label{eq:eq39}
\left|E_{t_n}^x[y_{t_{n+j}}\Delta W_{t_n,t_{n+j}}]-\hat E_{t_n}^x[y_{t_{n+j}}\Delta W_{t_n,t_{n+j}}]\right|\leq C(jh)^L.
\enq

Next let us consider the multi-dimensional case, i.e., $m>1, d>1$.
For a $d$-dimensional function $g:\R^d\to\R$, Gauss-Hermite quadratire formula for the integration $\int_{\R^d}{g(\x)e^{-\x^{\mathbf{T}}\x}d\x}$ is
\beq\label{eq:eq310} \int_{\R^d}{g(\x)e^{-\x^{\mathbf{T}}\x}d\x}\approx\sum_{i_1=1,\cdots,i_d=1}^{L,\cdots,L}{\omega_{i_1}\cdots\omega_{i_d}g(\mathbf{a}_\mathbf{i})}\enq
where $\mathbf{i}=(i_1,i_2,\cdots,i_d)$ and $\mathbf{a_i}=(a_{i_1},a_{i_2},\cdots,a_{i_d})$.
Similar to the way to deduce \cref{eq:eq36} and \cref{eq:eq39}, we obtain
\beq\label{eq:eq311}\left|E_{t_n}^x[y_{t_{n+j}}^{(l_1)}]-\frac{1}{\pi^{d/2}}\sum_{i_1=1,\cdots,i_d=1}^{L,\cdots,L}{\omega_{i_1}\cdots\omega_{i_d}y_{t_{n+j}}^{(l_1)}(\x+\sqrt{2jh}\mathbf{a_i})}\right|\leq C(jh)^L.\enq
\beq\label{eq:eq312}\left|E_{t_n}^x[y_{t_{n+j}}^{(l_1)}\Delta W_{t_n,t_{n+j}}^{(l_2)}]-\frac{\sqrt{2jh}}{\pi^{d/2}}\sum_{i_1=1,\cdots,i_d=1}^{L,\cdots,L}{\omega_{i_1}\cdots\omega_{i_d}y_{t_{n+j}}^{(l_1)}(\x+\sqrt{2jh}\mathbf{a_i})a_{l_2}}\right|\leq C(jh)^L.\enq
where $1\leq l_1 \leq m, 1\leq l_2 \leq d$ and 
\[
y_{t_{n+j}}=\left(y_{t_{n+j}}^{(1)},\cdots,y_{t_{n+j}}^{(m)}\right),
\Delta W_{t_{n+j}}=\left(\Delta W_{t_{n+j}}^{(1)},\cdots,\Delta W_{t_{n+j}}^{(d)}\right).
\]
\par From \cref{eq:eq310},  one can see that the computational cost increase exponentially in proportion to $L^d$, and obviously smaller $L$ is better. 
\par Based on the above discussions, we conclude that the error for the approximation of conditional expectation using Gauss-Hermite quadrature with $L$ sample points is $O(h^L)$. Therefore we state that it is enough to use $k$-point Gauss-Hermite quadrature when the local truncation error for time-discretization is $O(h^k)$.
\begin{remark}\label{remark31}
This idea can be applied to the other numerical schemes that use Gauss-Hermite quadrature for the approximation of conditional expectation. For instance, in the case of Crank-Nicolson scheme, the local truncation error for time-discretization is $O(h^3)$ and 3-point quadrature($L=3$) is enough. Note that in the literature $L=8$ is used for Crank-Nicolson scheme, this is redundant from the above discussion.
\end{remark}
\subsection{Construction of nested spatial grid}\label{sec:sec32}
In this subsection we discuss how to avoid spatial interpolation.
\par Let us assume that $m=d=1$ for the sake of simplicity. From \cref{eq:eq34} and \cref{eq:eq38}, we know that $\{y^{n+j}(x+\sqrt{2jh}a_i)\},(j=1\cdots k,i=1\cdots L)$ are needed to obtain $y^n(x)$ at time-space point $(t_n,x)(x\in D^n)$. In the case of $L>3$, the roots of Hermite polynomial of degree $L$, i.e., $\{a_i\}_{i=1}^L$ , are distributed nonlinearly and it seems to be impossible to construct the spatial grid $D^n$ so that it contains all the quadrature points needed. Therefore $D^n$ is usually constructed with equidistant points and $y^{n+j}(x+\sqrt{2jh}a_i)$ is approximated by $\hat y^{n+j}(x+\sqrt{2jh}a_i)$ where $\hat y^{n+j}(\cdot)$ is an interpolation based on $\{\left(x,y^{n+j}(x)\right)|x\in D^{n+j}\}$. Lagrange polynomial interpolation is often used for this purpose. The order of interpolation polynomial $r$ should be well-chosen because it influences the stability, accuracy and complexity of the scheme.(see \cite{Zhao14, Zhao17}) In \cite{Zhao14}, they chose small $r$(such as
$r = 4$ or $r = 6$) for lower-order schemes ($k\leq 3$) and bigger $r$ for higher-order schemes.(e.g.,
$r = 10$ or $r = 15$). We also note that the computational cost for the multi-dimensional interpolation is too expensive.
From the above discussion it is desirable not to use spatial interpolation.
\par For this purpose, we first construct the nested spatial grids by taking $L=3$.
In this case we have
\beq\label{eq:eq313} \{a_i|i=1,2,3\}=\{-\sqrt\frac{3}{2},0,\sqrt\frac{3}{2}\}, 
\{\omega_i|i=1,2,3\}=\{\frac{\sqrt\pi}{6},\frac{2\sqrt\pi}{3},\frac{\sqrt\pi}{6}\}.\enq
Now we set $\Delta x_0=\sqrt{2h}\cdot a_3=\sqrt{3h}$ and define $D^n$ by
\[D^n=\{x_l|x_l=l\Delta x_0, l=0,\pm 1,\cdots,\pm n\}.\]
Then for any $x_l\in D^n(l\in \mathbb{Z})$, we have $x_l+\sqrt{2h}a_i=x_{l+i-2}$ and 
\beq\label{eq:eq314} \{0\}=D^0\subset\cdots\subset D^n\subset D^{n+1} \subset \cdots\subset D^N.\enq
 The spatial grid satisfying \cref{eq:eq314} is called \emph{nested}.
\par Furthermore, we note that for any $j\in\mathbb{N}$ we have
\beq\label{eq:eq315}x_l\in D^n\Rightarrow x_l+j\sqrt{2h}a_i=x_{l+j(i-2)}\in D^{n+j}\subset D^{n+j^2}\enq
and $\#\{D^N\}=2N+1$, where $\#\{\cdot\}$ denotes a number of elements.
\par For the multi-dimensional case( i.e., $d>1$), if we define the spatial grid by
\beq\label{eq:eq316} D^n=\left\{x_\mathbf{l}|x_\mathbf{l}=(l_1\Delta x_0,\cdots,l_d\Delta x_0), |l_j|\leq n,l_j\in\mathbb{Z},j=1,\cdots,d\right\}\enq
then it is nested, i.e.,
\beq\label{eq:eq317} \{\mathbf{0}\}=D^0\subset\cdots\subset D^n\subset D^{n+1} \subset\cdots\subset D^N.\enq
Furthermore, for any $j\in\mathbb{N}$ we have
\beq\label{eq:eq318}\x_l\in D^n\Rightarrow \{\x+j\sqrt{2h}\mathbf{a_i}|i_j=1,2,3,j=1,\cdots,d\}\subset D^{n+j}\subset D^{n+j^2}\enq
and $\#\{D^N\}=(2N+1)^d$ where $\x+j\sqrt{2h}\mathbf{a_i}=(x_1+j\sqrt(2h)a_{i_1},\cdots,x_d+j\sqrt(2h)a_{i_d})$.
Using this nested spatial grid we propose a new kind of scheme that needs no spatial interpolation in the following.
\subsection{An efficient third-order scheme}\label{sec:sec33}
In this subsection we propose an efficient third-order scheme via the derivative approximation based on non-equidistant sample points.
\par Let $u(t):\R\to\R$ be $k+1$ times differentiable, then the Lagrange interpolation polynomial based on values $\{u_0,u_1,u_{2^2},\cdots,u_{k^2}\}$ on non-equidistant sample points $\{t_0,t_0+h,t_0+2^2h,\cdots,t_0+k^2h\}$ can be written as
\beq\label{eq:eq319} L(t)=\sum_{i=0}^k{\frac{\prod_{i\neq j}{(t-t_0-j^2h)}}{\prod_{i\neq j}{(i^2-j^2)}}h^{-k}u_{i^2}}\enq
and similar to the way we obtain \cref{eq:eq212}, we have
\begin{multline}\label{eq:eq320}
\begin{aligned}
\frac{dL}{dt}(t_0)&=\sum_{i=0}^k{\frac{(-1)^{k-1}h^{-1}u_{i^2}}{\prod_{j\neq i}{(i^2-j^2)}}\sum_{j\neq i}{\prod_{l\neq i,l\neq j}{l^2}}}\\&=h^{-1}\left(-\sum_{j\neq 0}{\frac{1}{j^2}u_0}+\sum_{i=1}^k{\frac{(-1)^{k-1}(k!)^2}{i^2\prod_{j\neq i}{(i^2-j^2)}}u_{i^2}}\right)
\end{aligned}
\end{multline}
Let $\beta_i^k$ be the coefficients of $u_{i^2}$ in \cref{eq:eq320}, then we have
\beq\label{eq:eq321} \frac{dL}{dt}(t_0)=\frac{\sum_{i=0}^k{\beta_i^k u_{i^2}}}{h}\enq
and
\beq\label{eq:eq322} \frac{du}{dt}(t_0)=\frac{\sum_{i=0}^k{\beta_i^k u_{i^2}}}{h}+O(h^k).\enq
We list $\{\beta_i^k\}$ for $k=1\cdots8$ in \Cref{tbl31} .
\begin{table}[tbhp]\label{tbl31}
\caption{$\{\beta_i^k\}$ in the derivative approximation based on the non-equidistant Lagrange interpolation}
\centering
\begin{tabular}{|c|c|c|c|c|c|c|c|c|c|} \hline
$\beta_i^k$& $i=0$ & $i=1$ & $i=2$ & $i=3$ & $i=4$ & $i=5$ & $i=6$ & $i=7$ & $i=8$ \\ \hline
$k=1$ &-1&1&&&&&&&\\ \hline
$k=2$ &-$\frac{5}{4}$&$\frac{4}{3}$&-$\frac{1}{12}$&&&&&&\\ \hline
$k=3$ &-$\frac{49}{36}$&$\frac{3}{2}$&-$\frac{3}{20}$&$\frac{1}{90}$ &&&&&\\ \hline
$k=4$ &-$\frac{205}{144}$&$\frac{8}{5}$&-$\frac{1}{5}$&$\frac{8}{315}$&-$\frac{1}{560}$&&&&\\ \hline
$k=5$ &-$\frac{5269}{3600}$&$\frac{5}{3}$&-$\frac{5}{21}$&$\frac{5}{126}$&-$\frac{5}{1008}$&$\frac{1}{3150}$&&&\\ \hline
$k=6$ &-$\frac{5369}{3600}$&$\frac{12}{7}$&-$\frac{15}{56}$&$\frac{10}{189}$&-$\frac{1}{112}$&$\frac{2}{1925}$&-$\frac{1}{16632}$&&\\ \hline
$k=7$&-$\frac{266681}{176400}$&$\frac{7}{4}$&-$\frac{7}{24}$&$\frac{7}{108}$&-$\frac{7}{528}$&$\frac{7}{3300}$&-$\frac{7}{30888}$&$\frac{1}{84084}$&\\ \hline
$k=8$&-$\frac{1077749}{705600}$&$\frac{16}{9}$&-$\frac{14}{45}$&$\frac{112}{1485}$&-$\frac{7}{396}$&$\frac{112}{32175}$&-$\frac{2}{3861}$&$\frac{16}{315315}$&-$\frac{1}{411840}$\\ \hline
\end{tabular}
\end{table}
\par From the stability theory of numerical ODEs, the stability of \Cref{scheme1} is related to the distribution of roots of the corresponding characteristic polynomial defined as following.(See \cite{Zhao14}) 
\beq\label{eq:eq323}P_\alpha^k(\lambda):=\sum_{i=0}^k{\alpha_i^k\lambda^{k-i}}.\enq
The roots $\{\lambda_i^k\}_{i=1}^k$ of this polynomial are said to satisfy the root conditions if $|\lambda_i^k|\leq 1$ and $|\lambda_i^k|=1\Rightarrow P_\alpha^k(\lambda_i^k)'\neq 0$.
In \cite{Zhao14}, they showed that the roots of $P_\alpha^k(\lambda)$ satisfy the root conditions for only $1\leq k\leq 6$.
\par If we approximate derivatives using \cref{eq:eq322}, the corresponding characteristic polynomial is 
\beq\label{eq:eq324}P_\beta^k(\lambda):=\sum_{i=0}^k{\beta_i^k\lambda^{k^2-i^2}}.\enq
We note that the roots $\{\gamma_i^k\}_{i=1}^k$ of \cref{eq:eq324} satisfy the root conditions for $1\leq k\leq 17$. (We obtained the table using Matlab 2013b which failed to calculate for $k>17$.) In \Cref{tbl32} we present the maximum absolute values of the roots except 1 for $k=2,\cdots,17$.

\begin{table}[tbhp]\label{tbl32}
\caption{The maximum absolute values of the roots of \cref{eq:eq324} except 1}
\centering
\begin{tabular}{|c|c|c|c|c|c|c|c|c|} \hline
$k$ & $ 2$ & $ 3$ & $ 4$ & $ 5$ & $ 6$ & $ 7$ & $ 8$ & $ 9$ \\ \hline
$max(|\lambda_{k,i}|)$&0.486&0.636&0.738&0.800&0.836&0.857&0.874&0.887\\ \hline
$k$ & $ 10$ & $ 11$ & $ 12$ & $ 13$ & $ 14$ & $ 15$ & $ 16$ & $ 17$ \\ \hline
$max(|\lambda_{k,i}|)$&0.896&0.903&0.901&0.91&0.921&0.926&0.931&0.935\\ \hline
\end{tabular}
\end{table}

\par
\par Now approximating the derivatives in \cref{eq:eq27} and \cref{eq:eq28} by \cref{eq:eq322}, we obtain another discrete scheme for BSDE \cref{eq:eq21} as following:
\begin{scheme}\label{scheme3}
Assume that $y^{N-i},z^{N-i}(i=0,\cdots,k^2-1)$ are known. For $n=N-k^2,\cdots,0$ solve $y^n=y^n(\mathbf{x}),z^n=z^n(\mathbf{x})$ at time-space point $(t_n,\mathbf{x})$ by
\beq\label{eq:eq325} z^n=\frac{\sum_{j=1}^k{\beta_j^kE_{t_n}^\mathbf{x}[y^{n+j^2}\Delta W_{t_n,t_{n+j^2}}^\mathbf{T}]}}{h}\enq
\beq\label{eq:eq326} -\beta_0^ky^n=\sum_{j=1}^k{\alpha_j^kE_{t_n}^\mathbf{x}[y^{n+j^2}]}+h\cdot f(t_n,y^n,z^n)\enq
\end{scheme}
The truncation error is $O(h^k)$ again.
Now we introduce the spatial grid $D^n$ for each $t_n$ as in \cref{eq:eq316} and obtain the following fully-discrete scheme. Note that we set the step size $k=3$ to balance the time discretization error and error from the approximation of conditional expectation.

\begin{scheme}\label{scheme4}
Assume that $y^{N-i},z^{N-i}$ on $D^{N-i}$ are known for $i=0,\cdots,8$. For $n=N-9,\cdots,0,\mathbf{x}\in D^n$ solve $y^n=y^n(\mathbf{x}),z^n=z^n(\mathbf{x})$ by
\beq\label{eq:eq327} z^n=\frac{\sum_{j=1}^3{\beta_j^3\hat{E}_{t_n}^\mathbf{x}[y^{n+j^2}\Delta W_{t_n,t_{n+j^2}}^\mathbf{T}]}}{h}\enq
\beq\label{eq:eq328} -\beta_0^3y^n=\sum_{j=1}^3{\beta_j^3\hat{E}_{t_n}^\mathbf{x}[y^{n+j^2}]}+h\cdot f(t_n,y^n,z^n)\enq
where 
\beq\label{eq:eq329}\hat{E}_{t_n}^\mathbf{x}[y^{n+j^2}] =\frac{1}{\pi^{d/2}}\sum_{i_1=1,\cdots,i_d=1}^{3,\cdots,3}{\omega_{i_1}\cdots\omega_{i_d}y^{n+j^2}(\x+j\sqrt{2h}\mathbf{a_i})}\enq
\beq\label{eq:eq330}\hat{E}_{t_n}^\mathbf{x}[y^{n+j^2}\Delta W_{t_n,t_{n+j^2}}^\mathbf{T}] =\frac{j\sqrt{2h}}{\pi^{d/2}}\sum_{i_1=1,\cdots,i_d=1}^{3,\cdots,3}{\omega_{i_1}\cdots\omega_{i_d}y^{n+j^2}(\x+j\sqrt{2h}\mathbf{a_i})\mathbf{a_i^T}}.\enq
\end{scheme}
Note that in \cref{eq:eq329} and \cref{eq:eq330}, $\x+j\sqrt{2h}\mathbf{a_i}\in D^{n+j}\subset D^{n+j^2}$ and no interpolation is needed.
\par Especially for the one dimensional case($m=d=1$) we have the following scheme.
\begin{scheme}\label{scheme5}
Assume that $y^{N-i},z^{N-i}$ on $D^{N-i}$ are known for $i=0,\cdots,8$. For $n=N-9,\cdots,0,l=0,\pm1,\cdots,\pm n, $ solve $y^{n,l}=y^n(x_l),z^{n,l}=z^n(x_l) (x_l=l\Delta x_0)$ by
\beq\label{eq:eq331} z^{n,l}=\sqrt{\frac{2}{\pi h}}\sum_{j=1}^3{j\cdot\beta_j^3\sum_{i=1}^3{\omega_i a_i y^{n+j^2,l+j(i-2)}}}\enq
\beq\label{eq:eq332} -\alpha_0^ky^{n,l}=\frac{1}{\sqrt\pi}\sum_{j=1}^3{\beta_j^3\sum_{i=1}^3{\omega_i y^{n+j^2,l+j(i-2)}}}+h\cdot f(t_n,y^{n,l},z^{n,l})\enq
where $a_i,\omega_i(i=1,2,3)$ are given by \cref{eq:eq313}.\\ 
\end{scheme}

\begin{remark}\label{remark32}
Unlike the \Cref{scheme1} and \Cref{scheme2}, \Cref{scheme3} needs $k^2$ initial approximations.
\end{remark}
\begin{remark}\label{remark33}
Note that in the above schemes procedures for solving $y^n(\x)$ at different grid points $\x\in D^n$ are independent and that parallel computing technique can be used. for the multi-dimensional case ($d>1$) sparse grid quadrature would be effective.(See \cite{Zhao10})
\end{remark}

\section{Numerical Experiments}\label{sec:sec4}
In this section we demonstrate the convergence and efficiency of the new scheme through some numerical examples. We obtained the numerical results with \texttt{Matlab 2013b} on a computer with \texttt{Intel(R) Core(TM) i7-3770 CPU @ 3.40GHz (8 CPUs)}. 
For all the examples we set the terminal time $T=1$ and measured the error at $t=0$. We used the true solution for initial approximations but it could be calculated numerically using other schemes such as Crank-Nicolson scheme.

\begin{example}\label{ex41}\end{example}
Let us consider the following one dimensional backward stochastic differential equation. The equation is from \cite{Zhao10}.
\beq\label{eq:eq41}
\left\{\begin{array}{l}
-dy_t=\frac{1}{2}\left[e^{t^2}-4ty_t-3e^{t^2-y_te^{-t^2}}+z_t^2e^{-t^2}\right]dt-z_tdW_t\\
y_T=\textnormal{ln}(\textnormal{sin}W_T+3)e^{T^2}
\end{array}\right.
\enq
The analytic solution of this equation is $(y_t,z_t)=\left(\textnormal{ln}(\textnormal{sin}W_t+3)e^{t^2},\frac{\textnormal{cos}{W_t}}{\textnormal{sin}{W_t}+3}e^{t^2}\right)$ and $(y_0,z_0)=(\textnormal{ln}{3},\frac{1}{3})$.
\par First we test the influence of the number of sample points of Gauss-Hermite quadrature on the error of scheme. To this end we investigated the convergence rate for different step size $k$ and number of quadrature points $L$. We calculated the solution by \Cref{scheme2} increasing $N$ from $2^2$ to $2^6$. When using polynomial interpolation of degree $r$, we set the diameter of spatial partition $\Delta x=(\Delta t)^{\frac{k+1}{r+1}}$ to balance the error from time discretization and that from spatial interpolation.(see \cite{Zhao10} for details) 
\begin{table}[tbhp]\label{tbl41}
\caption{Convergence rate and running time of \Cref{scheme2} for \Cref{ex41}, for different number of quadrature points($L$)}
\centering
\begin{small}
\begin{tabular}{|c|c|c|c|c|c|c|} \hline
\multicolumn{2}{|c|}{}& $k=2$ & $k=3$ & $k=4$ & $k=5$ & $k=6$\\
\multicolumn{2}{|c|}{\Cref{scheme2}}& $r=5$ & $r=8$ & $r=16$ & $r=18$ & $r=20$ \\ \hline
& $\textnormal{CR}_y$ &1.826&-0.406&NaN&NaN&NaN\\ 
$L=2$ & $\textnormal{CR}_z$ &2.148&-0.034&NaN&NaN&NaN\\ 
&RT&4.57s&11.62s&61.70s&95.83s&143.08s\\ \hline
& $\textnormal{CR}_y$ &1.835&2.597&3.186&4.556&3.307\\ 
$L=3$ & $\textnormal{CR}_z$ &2.087&3.110&4.173&5.025&1.004\\ 
&RT&6.36s&20.40s&90.12s&145.24s&213.92s\\ \hline
& $\textnormal{CR}_y$ &1.835&2.612&3.173&NaN&NaN\\ 
$L=4$ & $\textnormal{CR}_z$ &2.108&4.104&5.131&NaN&NaN\\ 
&RT&10.15s&25.52s&122.89s&207.03s&308.21s\\ \hline
& $\textnormal{CR}_y$ &1.835&2.613&3.186&4.556&4.208\\ 
$L=5$ & $\textnormal{CR}_z$ &2.116&4.225&4.220&4.319&3.095\\ 
&RT&12.21s&36.39s&164.56s&259.48s&389.90s\\ \hline
& $\textnormal{CR}_y$ &1.835&2.613&3.185&4.554&5.377\\ 
$L=6$ & $\textnormal{CR}_z$ &2.114&4.216&4.200&4.713&5.620\\ 
&RT&17.36s&42.79s&194.90s&311.60s&466.25s\\ \hline
\end{tabular}
\end{small}
\end{table}
\begin{table}[tbhp]\label{tbl42}
\caption{Comparison of \Cref{scheme2} and \Cref{scheme5} for \Cref{ex41}}
\centering
\begin{small}
\begin{tabular}{|c|c|c|c|c|c|c|c|} \hline
\multicolumn{2}{|c|}{}& $N=2^4$ & $N=2^5$ & $N=2^6$ & $N=2^7$ & $N=2^8$ & CR\\ \hline
& $Err_y$ & 2.486E-03&3.440E-04&4.507E-05&5.761E-06&7.282E-07&2.94\\ 
\Cref{scheme2} & $Err_z$ & 4.779E-05&1.696E-07&4.013E-07&7.878E-08&1.170E-08&2.51\\ 
&RT&2.61s&12.27s&50.88s&211.64s&874.89s&\\ \hline
& $Err_y$ &1.064E-04&2.368E-05&3.670E-06&5.022E-07&6.516E-08&2.69\\ 
\Cref{scheme5} & $Err_z$ &9.661E-05&1.846E-05&2.728E-06&3.685E-07&4.783E-08&2.76\\ 
&RT&0.02s&0.05s&0.27s&1.40s&6.98s&\\ \hline
\end{tabular}
\end{small}
\end{table}
\par We choose the degree of interpolation polynomial $r$ carefully for each $k$ so that the scheme converges in a stable manner. In \Cref{tbl41}, we listed the convergence rate for $y$ and $z$(namely $\textnormal{CR}_y$ and $\textnormal{CR}_z$) and the running time (`RT' in the table) in seconds. We calculate the convergence rate using the least square linear fitting. From the \cref{tbl41}, we see that $k$-points Gauss-Hermite quadrature is enough for $k$-step scheme where the local time discretization error is $O(h^k)$.
\par Next we compare the error, the convergence rate and the running time of \Cref{scheme2}\\(with $k=3,r=8,L=8$) and \Cref{scheme5} for \cref{eq:eq41}. We list the experiment result in \Cref{tbl42}. It shows that \Cref{scheme5} is much faster than \Cref{scheme2} while achieving the same convergence rate.

\begin{example}\label{ex42}\end{example}
Here we test the following multi-dimensional backward stochastic differential equation\\($m>1,d=1$) 
\beq\label{eq:eq42}
\left\{\begin{array}{l}
-dy_t=\mathbf{A}y_t\Vert y_t\Vert^2dt-z_tdW_t\\
y_T=\left(\textnormal{sin}(W_T+T),\textnormal{cos}(W_T+T)\right)^\textnormal{T}
\end{array}\right.
\enq
where $y_t=(y_t^1,y_t^2)^\textnormal{T}, \mathbf{A}=\begin{pmatrix}0.5&-1\\1&0.5\end{pmatrix}, \Vert y_t \Vert ^2=(y_t^1)^2+(y_t^2)^2$.(This is also from \cite{Zhao10})
\par The analytic solution of this equation is 
\[ y_t=\begin{pmatrix}\textnormal{sin}(W_t+t)\\\textnormal{cos}(W_t+t)\end{pmatrix}, 
 z_t=\begin{pmatrix}\textnormal{cos}(W_t+t)\\-\textnormal{sin}(W_t+t)\end{pmatrix}, \]
 and $y_0=\begin{pmatrix}0\\1\end{pmatrix},z_0=\begin{pmatrix}1\\0\end{pmatrix}$.
As in the above example we test the \Cref{scheme2} and \Cref{scheme4} and list the result in \cref{tbl43}.

\begin{table}[tbhp]\label{tbl43}
\caption{Comparison of \Cref{scheme2} and \Cref{scheme5} for \Cref{ex42}}
\centering
\begin{small}
\begin{tabular}{|c|c|c|c|c|c|c|c|} \hline
\multicolumn{2}{|c|}{}& $N=2^4$ & $N=2^5$ & $N=2^6$ & $N=2^7$ & $N=2^8$ & CR\\ \hline
& $Err_y$ &1.606E-04&2.199E-05&2.910E-06&3.719E-07&4.694E-08&2.94\\ 
\Cref{scheme2} & $Err_z$ &6.614E-04&8.724E-05&1.122E-05&1.430E-06&1.793E-07&2.96\\ 
&RT&5.19s&24.94s&105.73s&455.76s&1908.70s&\\ \hline
& $Err_y$ &6.108E-04&1.401E-04&2.204E-05&3.058E-06&4.019E-07&2.79\\ 
\Cref{scheme5} & $Err_z$ &1.915E-03&2.176E-04&2.592E-05&3.197E-06&3.989E-07&3.03\\ 
&RT&0.01s&0.11s&0.37s&1.69s&8.70s&\\ \hline
\end{tabular}
\end{small}
\end{table}

\begin{example}\label{ex43}\end{example}
In this example, we consider the case of $m=1,d>1$. Let $W_t=\begin{pmatrix}W_t^1\\W_t^2\end{pmatrix}$ be a two-dimensional Brownian motion where $W_t^1$ and $W_t^2$ are independent and standard one-dimensional Brownian motions.
\beq\label{eq:eq43}
\left\{\begin{array}{l}
-dy_t=(y_t-z_t^1-2tz_t^2)dt-z_tdW_t\\
y_T=\textnormal{sin}(W_T^1+T^2)\textnormal{cos}(W_T^2+T)
\end{array}\right.
\enq
where $z_t=(z_t^1,z_t^2)$.
\par The analytic solution of is 
\[ y_t=\textnormal{sin}(W_t^1+t^2)\textnormal{cos}(W_t^2+t), 
 z_t=\begin{pmatrix}\textnormal{cos}(W_t^1+t^2)\textnormal{cos}(W_t^2+t)\\-\textnormal{sin}(W_t^1+t^2)\textnormal{sin}(W_t^2+t)\end{pmatrix}^\textnormal{T}\]
 and $y_0=0,z_0=(1,0)$.
Again we test the \Cref{scheme2} and \Cref{scheme4} solving \cref{eq:eq43} and list the result in \Cref{tbl44}. One can see that the \Cref{scheme2} fails to converge and its running time is very long. We believe that this is because of the multi-dimensional interpolation. 
\\
\par We also note that our new scheme is of third-order for all examples and very efficient especially for the multi-dimensional problems.
\begin{table}[tbhp]\label{tbl44}
\caption{Comparison of \Cref{scheme2} and \Cref{scheme5} for \Cref{ex43}}
\centering
\begin{small}
\begin{tabular}{|c|c|c|c|c|c|c|} \hline
\multicolumn{2}{|c|}{}& $N=2^4$ & $N=2^5$ & $N=2^6$ & $N=2^7$ & CR\\ \hline
& $Err_y$ &9.597e-02&4.198e-02&5.409e+05&-& -11.21\\ 
\Cref{scheme2} & $Err_z$ &1.929e-01&5.463e-02&1.703e+06&-&-11.53\\ 
&RT&139.33s&3252.48s&67227.23s&-&\\ \hline
& $Err_y$ &2.464E-04&1.010E-04&1.332E-05&1.326E-06&2.55\\ 
\Cref{scheme5} & $Err_z$ &4.359E-03&1.578E-03&3.024E-04&4.617E-05&2.21\\ 
&RT&0.32s&3.99s&51.56s&851.71s&\\ \hline
\end{tabular}
\end{small}
\end{table}

\section{Discussion and Conclusions}\label{sec:sec5}
In this paper we proposed a new efficient third-order numerical scheme for BSDEs. Our main contributions are related to the approximation of conditional expectation using Gauss-Hermite quadrature. First we carried out analysis on the number of quadrature points needed. We saw that $k$ quadrature points are enough for $k$-step scheme where the time discretization error is $O(h^k)$. Next we proposed a new kind of efficient scheme. Our scheme is based on the nested spatial grid and the approximation of derivative using non-equidistant sample points. The proposed scheme does not include the spatial interpolation, which is very costly in solving multi-dimensional problems numerically. As mentioned in \Cref{remark31}, some ideas can be extended to other types of schemes, e.g., Crank-Nicolson scheme and so on. We note that our scheme needs more initial data than the former schemes. The numerical examples show that our scheme is of third order and very efficient, especially for the multi-dimensional problems.

\bibliographystyle{siamplain}
\bibliography{}

\end{document}